\documentclass{icm2010}

\title[Rigidity for von Neumann algebras and their invariants]{Rigidity for von Neumann algebras and their invariants}
\author[Stefaan Vaes]{Stefaan Vaes
\thanks{Partially
    supported by ERC Starting Grant VNALG-200749, Research
    Programme G.0231.07 of the Research Foundation --
    Flanders (FWO) and K.U.Leuven BOF research grant OT/08/032.}}
\contact[stefaan.vaes@wis.kuleuven.be]{K.U.Leuven\\ Department of Mathematics\\ Celestijnenlaan 200B\\ B-3001 Leuven (Belgium)}

\newtheorem{theorem}{Theorem}[section]

\newtheorem{lemma}[theorem]{Lemma}

\theoremstyle{definition}
\newtheorem{definition}[theorem]{Definition}
\newtheorem{example}[theorem]{Example}
\newtheorem{problem}{Problem}

\renewcommand{\thesubsection}{{\normalsize \arabic{section}.\arabic{subsection}}}
\newcommand{\smalltitle}[1]{\subsection{\normalsize #1}}
\newcommand{\smalltitlezonder}[1]{\subsection*{\normalsize #1}}

\newcommand{\F}{\mathbb{F}}
\newcommand{\R}{\mathbb{R}}
\newcommand{\SL}{\operatorname{SL}}
\newcommand{\actson}{\curvearrowright}
\newcommand{\rL}{\mathord{\text{\rm L}}}
\newcommand{\B}{\mathord{\text{\rm B}}}
\newcommand{\ovt}{\overline{\otimes}}
\newcommand{\cH}{\mathcal{H}}
\newcommand{\ot}{\otimes}
\newcommand{\recht}{\rightarrow}
\newcommand{\C}{\mathbb{C}}
\newcommand{\cF}{\mathcal{F}}
\newcommand{\Z}{\mathbb{Z}}
\newcommand{\T}{\mathbb{T}}
\newcommand{\PSL}{\operatorname{PSL}}
\newcommand{\Aut}{\operatorname{Aut}}
\newcommand{\Inn}{\operatorname{Inn}}
\newcommand{\Ad}{\operatorname{Ad}}
\newcommand{\Out}{\operatorname{Out}}
\newcommand{\cU}{\mathcal{U}}
\newcommand{\bim}[3]{\mathord{\raisebox{-0.4ex}[0ex][0ex]{\scriptsize $#1$}{#2}\hspace{-0.2ex}\raisebox{-0.4ex}[0ex][0ex]{\scriptsize $#3$}}}
\newcommand{\Bimod}{\operatorname{Bimod}}
\newcommand{\vphi}{\varphi}
\newcommand{\M}{\operatorname{M}}
\newcommand{\id}{\mathord{\text{\rm id}}}
\newcommand{\al}{\alpha}
\newcommand{\N}{\mathbb{N}}
\newcommand{\otalg}{\otimes_{\text{\rm alg}}}
\newcommand{\cK}{\mathcal{K}}
\newcommand{\eps}{\varepsilon}
\newcommand{\lspan}{\operatorname{span}}
\newcommand{\invlimit}{\varprojlim}
\newcommand{\cbnorm}[1]{\| #1 \|_{\text{\rm cb}}}
\newcommand{\cG}{\mathcal{G}}
\newcommand{\om}{\omega}
\newcommand{\Tr}{\operatorname{Tr}}
\newcommand{\cN}{\mathcal{N}}
\newcommand{\cR}{\mathcal{R}}
\newcommand{\module}{\operatorname{mod}}
\newcommand{\Centr}{\operatorname{Centr}}
\newcommand{\nutil}{\widetilde{\nu}}

\newcommand{\Scentr}{\mathcal{S}_{\text{\rm centr}}}
\newcommand{\Rep}{\operatorname{Rep}}

\begin{document}

\begin{abstract}
We give a survey of recent classification results for von Neumann algebras $\rL^\infty(X) \rtimes \Gamma$ arising from measure preserving group actions on probability spaces. This includes II$_1$ factors with uncountable fundamental groups and the construction of W$^*$-superrigid actions where $\rL^\infty(X) \rtimes \Gamma$ entirely remembers the initial group action $\Gamma \actson X$.
\end{abstract}



\maketitle

\section{Classifying II$_1$ factors, a panoramic overview}\label{sec.panorama}

A von Neumann algebra is a an algebra of bounded linear operators on a Hilbert space that is closed under the adjoint $*$-operation and that is closed in the weak operator topology. Von Neumann algebras arise naturally in the study of groups and their actions on measure spaces. These constructions go back to Murray and von Neumann's seminal papers \cite[Chapter XII]{MvN36} and \cite[\S 5.3]{MvN43}.

\begin{itemize}
\item If $\Gamma$ is a countable group, the left translation unitary operators on $\ell^2(\Gamma)$ generate the \emph{group von Neumann algebra} $\rL \Gamma$.
\item Every action $\Gamma \actson (X,\mu)$ of a countable group $\Gamma$ by measurable transformations of a measure space $(X,\mu)$ and preserving sets of measure zero, gives rise to the \emph{group measure space von Neumann algebra} $\rL^\infty(X) \rtimes \Gamma$.
\end{itemize}

It is a central problem in the theory of von Neumann algebras to classify $\rL \Gamma$ and $\rL^\infty(X) \rtimes \Gamma$ in terms of the group $\Gamma$ or the group action $\Gamma \actson (X,\mu)$. More generally, classifying or distinguishing families of von Neumann algebras is extremely challenging. In the first part of this exposition, I give a panoramic overview of the spectacular progress that has been made in this area over the last years. The overview is more thematically ordered than chronologically and necessarily incomplete. Several related important topics, including Jones' theory of subfactors or Voiculescu's free probability theory, are not treated.

\vspace{0.5ex}

\noindent{\it All notions that are written in italics are defined in the preliminary section \ref{sec.prelim}.}

\smalltitle{II$_1$ factors}

The \lq simple\rq\ von Neumann algebras $M$ are those that cannot be written as a direct sum of two. Equivalently, the center of $M$ is trivial and $M$ is called a factor. Murray and von Neumann have classified factors into three types \cite{MvN36} and proven that
every von Neumann algebra can be decomposed as a direct integral of factors \cite{vN49}. Connes \cite{Co72} showed how general factors can be built up from those that admit a finite positive \emph{trace,} called \emph{II$_1$ factors.} The final form of this decomposition theory is due to Connes and Takesaki \cite{Ta73,CT76}. Altogether, II$_1$ factors form the basic building blocks of arbitrary von Neumann algebras.

The \emph{group von Neumann algebras} $\rL \Gamma$ always admit a finite positive trace and are factorial if and only if $\Gamma$ has infinite conjugacy classes (icc). When $\Gamma \actson (X,\mu)$ is essentially free, ergodic\footnote{Essential freeness means that for every $g \neq e$, the set of $x \in X$ with $g \cdot x = x$ has measure zero. Ergodicity means that globally $\Gamma$-invariant measurable subsets have either measure $0$ or a complement of measure $0$.} and probability measure preserving (p.m.p.), the \emph{group measure space von Neumann algebra} $\rL^\infty(X) \rtimes \Gamma$ is a II$_1$ factor. Moreover, as proven by Singer \cite{Si55}, its isomorphism class only depends on the \emph{equivalence relation} given by the orbits of $\Gamma \actson (X,\mu)$. This lead to the study of group actions up to orbit equivalence \cite{Dy58} and we refer to \cite{Sh05,Po06b,Fu09,Ga10} for surveys of the recent developments in this area of ergodic theory.

\smalltitle{(Non)-isomorphism of II$_1$ factors}\label{par.noniso}

Two von Neumann algebras can be isomorphic in unexpected ways.  While all hyperfinite\footnote{A von Neumann algebra is hyperfinite if it is the direct limit of finite dimensional subalgebras.} II$_1$ factors were already shown to be isomorphic in \cite{MvN43}, the culmination came with Connes' uniqueness theorem for \emph{amenable} II$_1$ factors \cite{Co75b} implying that all $\rL \Gamma$ and all $\rL^\infty(X) \rtimes \Gamma$ are isomorphic when $\Gamma$ is any amenable icc group or $\Gamma \actson (X,\mu)$ is an arbitrary free ergodic p.m.p.\ action of an amenable group.

In the early years examples of non-isomorphic II$_1$ factors $M$ were obtained by analyzing asymptotically central sequences\footnote{A bounded sequence $(x_n)$ in a II$_1$ factor $M$ is called asymptotically central if $x_n y - y x_n$ converges to $0$ in the strong operator topology for every $y \in M$.} of elements in $M$~: property Gamma\footnote{A II$_1$ factor has property Gamma if it admits an asymptotically central sequence of unitaries having trace $0$.} \cite{MvN43} allowed to prove that the free group factors $\rL \F_n$ are not hyperfinite, a refinement yielded uncountably many non-isomorphic II$_1$ factors \cite{McD69} and the $\chi$-invariant \cite{Co75a} provided the first examples where $M$ is non-isomorphic to its opposite algebra $M^{\text{\rm op}}$ and where $M \not\cong M \ovt M$.

The first rigidity phenomena for von Neumann algebras were discovered by Connes \cite{Co80} who showed that the \emph{fundamental group\footnote{The fundamental group $\cF(M)$ of a II$_1$ factor $M$ consists of the numbers $\tau(p)/\tau(q)$ where $p$ and $q$ run over the projections in $M$ satisfying $pMp \cong qMq$. Here $\tau$ denotes the trace on $M$.}} of $\rL \Gamma$ is countable when $\Gamma$ is an icc property (T) group. Several properties of groups -- including \emph{property (T),} the Haagerup property and related approximation properties -- were shown in \cite{CJ83,CH88,Jo00} to actually be properties of $\rL \Gamma$, leading to remarkable non-isomorphism and non-embeddability theorems for group von Neumann algebras. Altogether it became clear that the world of II$_1$ factors is extremely rich, but that understanding the natural examples $\rL \Gamma$ or $\rL^\infty(X) \rtimes \Gamma$ in terms of the initial group or action is intrinsically very difficult.

\smalltitle{Popa's deformation/rigidity theory}\label{par.defrig}

A major breakthrough in the classification of II$_1$ factors was realized by Popa and his discovery of deformation/ rigidity theory \cite{Po01} (see \cite{Po06b,Va06a} for a survey). Typically, Popa studies von Neumann algebras $M$ that have a rigid subalgebra -- e.g.\ given by the \emph{relative property (T)} -- such that the \lq complement\rq\ has a strong deformation property. This gives the rigid subalgebra a canonical position within the ambient von Neumann algebra and has lead in \cite{Po01} to the first example of a II$_1$ factor with trivial \emph{fundamental group:} $M=\rL(\Z^2 \rtimes \SL(2,\Z))= \rL^\infty(\T^2) \rtimes \SL(2,\Z)$. The canonical position of $\rL^\infty(\T^2)$ implies that the fundamental group of $M$ equals the fundamental group of the \emph{orbit equivalence relation} of the action $\SL(2,\Z) \actson \T^2$, which is trivial because of \cite{Ga99,Ga01}.

In \cite{Po03,Po04} Popa established a striking progress in his deformation/rigidity program by proving the following strong rigidity theorem for group measure space factors. Take an arbitrary free ergodic p.m.p.\ action $\Gamma \actson (X,\mu)$ of a property (T) group $\Gamma$ and let $\Lambda \actson (Y,\eta) := (Y_0,\eta_0)^\Lambda$ be the \emph{Bernoulli action} of an arbitrary icc group $\Lambda$. If the group measure space factors $\rL^\infty(X) \rtimes \Gamma$ and $\rL^\infty(Y) \rtimes \Lambda$ are isomorphic, then $\Gamma$ must be isomorphic with $\Lambda$ and their actions must be conjugate. Popa's strong rigidity theorem was the first result ever where conjugacy of actions could be deduced from the mere isomorphism of group measure space factors.

\smalltitle{Fundamental groups of II$_1$ factors}\label{par.fundamental}

Progress in the classification of group measure space factors went hand in hand with major developments in the calculation of invariants of II$_1$ factors. The most well known invariant of $M$ is the fundamental group $\cF(M)$ of Murray and von Neumann \cite{MvN43}. In one of their long-standing questions they asked what subgroups of $\R^*_+$ might occur as $\cF(M)$.

Until 10 years ago, progress on this question has been scarce. Some II$_1$ factors, including the hyperfinite II$_1$ factor \cite{MvN43} and $\rL \F_\infty$ \cite{Vo89,Ra91}, were shown to have fundamental group $\R^*_+$ while Connes \cite{Co80} proved that $\cF(\rL \Gamma)$ is countable whenever $\Gamma$ is an icc property (T) group. A breakthrough in the understanding of fundamental groups came with Popa's first examples of II$_1$ factors having trivial fundamental group \cite{Po01} and having prescribed countable fundamental group \cite{Po03}. It remained a major open problem whether uncountable groups $\neq \R^*_+$ could appear as fundamental group.

In \cite{PV08a} we solved this problem and proved that $\cF(\rL^\infty(X) \rtimes \F_\infty)$ ranges over a large family of subgroups of $\R^*_+$, including all countable subgroups and many uncountable subgroups that can have arbitrary Hausdorff dimension between $0$ and $1$. A similar result is true \cite{PV08c} when $\F_\infty$ is replaced by almost any infinite free product of non-trivial groups, while $\cF(\rL^\infty(X) \rtimes \Gamma)$ is necessarily trivial when $\Gamma \actson (X,\mu)$ is an arbitrary free ergodic p.m.p.\ action of a free product of two finitely generated groups, one of them having property (T).
So far, \cite{PV08a,PV08c} provide the only known constructions of group measure space factors with fundamental group different from $\{1\}$ or $\R^*_+$.

The fundamental group $\cF(M)$ of a II$_1$ factor $M$ can also be viewed as the set of $t > 0$ such that the II$_\infty$ factor $M \ovt \B(H)$ admits an automorphism scaling the (infinite) trace $\tau \ot \Tr$ by $t$. In \cite{PV08a} we provide examples where $\cF(M) = \R^*_+$, although $M \ovt \B(H)$ admits no continuous trace-scaling action of $\R^*_+$.

\smalltitle{Outer automorphisms and generalized symmetries}

Another invariant of a II$_1$ factor $M$ is its outer automorphism group\footnote{The outer automorphism group $\Out M$ is defined as the quotient $\Aut M /\Inn M$, where $\Inn M$ denotes the normal subgroup of $\Aut M$ consisting of the inner automorphisms $\Ad u$, $u \in \cU(M)$.} $\Out M$.
In \cite{IPP05} Ioana, Peterson and Popa established Bass-Serre isomorphism and subgroup (rather subalgebra) theorems for amalgamated free products of von Neumann algebras. As a consequence they obtained the first calculations of outer automorphism groups and proved that $\Out M$ can be any compact abelian group. In particular, they positively answered the question on the existence of II$_1$ factors without outer automorphisms. Later we showed \cite{FV07} that in fact $\Out M$ can be any compact group. The results in \cite{IPP05,FV07} are existence theorems involving a Baire category argument. We obtained the first concrete and explicit calculations of $\Out M$ in \cite{PV06,Va07} and proved that $\Out M$ can be any countable group.

Both the elements of the fundamental group and the automorphisms of a II$_1$ factor $M$ give rise to \emph{Hilbert $M$-$M$-bimodules} $\bim{M}{\cH}{M}$. An $M$-$M$-bimodule that is finitely generated, both as a left and as a right $M$-module is said to be of finite Jones index \cite{Jo82}. The finite index $M$-$M$-bimodules form a C$^*$-tensor category $\Bimod M$ and this should be considered as the generalized (or quantum) symmetry group of $M$. Whenever $M \subset P$ is a finite index subfactor, $\bim{M}{\rL^2(P)}{M}$ is a finite index bimodule. In this sense, $\Bimod M$ also encodes the subfactor structure of $M$. In \cite{Va06b} I proved the existence of II$_1$ factors $M$ such that $\Bimod M$ is trivial, i.e.\ only consists of multiples of the trivial bimodule $\rL^2(M)$. Such II$_1$ factors have trivial fundamental group, trivial outer automorphism group and no non-trivial finite index subfactors. Explicit examples were provided in \cite{Va07} where also several concrete calculations of $\Bimod M$ were made. These calculations were exploited in \cite{DV10} to give a full classification of all finite index subfactors of certain II$_1$ factors. In \cite{FV08} every representation category of a compact group $K$ is realized as $\Bimod M$. More precisely, for every compact group $K$ we prove the existence of a minimal action of $K$ on a II$_1$ factor $M$ such that, denoting by $M^K$ the subfactor of $K$-invariant elements, the natural faithful tensor functor $\Rep K \recht \Bimod(M^K)$ is \lq surjective\rq, i.e.\ an equivalence of categories.

\smalltitle{W$^*$-superrigidity and uniqueness of Cartan subalgebras}\label{par.superrigid}

Two free ergodic p.m.p.\ actions $\Gamma \actson (X,\mu)$ and $\Lambda \actson (Y,\eta)$ are called
\begin{itemize}
\item conjugate (or isomorphic), if there exists an isomorphism of probability spaces $\Delta : X \recht Y$ and an isomorphism of groups $\delta : \Gamma \recht \Lambda$ satisfying $\Delta(g \cdot x) = \delta(g) \cdot \Delta(x)$ for all $g \in \Gamma$ and a.e.\ $x \in X$;
\item orbit equivalent (OE), if there exists an isomorphism of probability spaces $\Delta : X \recht Y$ satisfying $\Delta(\Gamma \cdot x) = \Lambda \cdot \Delta(x)$ for a.e.\ $x \in X$;
\item W$^*$-equivalent, if $\rL^\infty(X) \rtimes \Gamma \cong \rL^\infty(Y) \rtimes \Lambda$.
\end{itemize}
Obviously, conjugacy of actions implies orbit equivalence and Singer \cite{Si55} proved that an orbit equivalence is the same as a W$^*$-equivalence sending the \emph{group measure space Cartan subalgebras\footnote{In general, a Cartan subalgebra $A$ of a II$_1$ factor $M$ is a maximal abelian subalgebra whose normalizing unitaries generate $M$. Whenever $\Gamma \actson (X,\mu)$ is a free ergodic p.m.p.\ action, $\rL^\infty(X) \subset \rL^\infty(X) \rtimes \Gamma$ is an example of a Cartan subalgebra that we call of group measure space type. Not all II$_1$ equivalence relations can be implemented by a free action of a countable group \cite{Fu98b} and hence a general Cartan subalgebra need not be of group measure space type.}} $\rL^\infty(X)$ and $\rL^\infty(Y)$ onto each other. Rigidity theory for group actions aims at establishing the converse implications under appropriate assumptions. Pioneering OE rigidity results were obtained by Zimmer \cite{Zi79,Zi84} and the first breakthrough W$^*$-rigidity theorems were proven by Popa \cite{Po01,Po03,Po04} (see paragraph \ref{par.defrig}). Our aim here however is to discuss the ideal kind of rigidity, labeled W$^*$- (respectively OE-) superrigidity, where the entire isomorphism class of $\Gamma \actson (X,\mu)$ is recovered from its W$^*$-class (resp.\ OE class).

While a striking number of OE superrigid actions have been discovered over the last 10 years \cite{Fu98b,Po05,Po06a,Ki06,Io08,PV08b,Ki09,PS09},
the first W$^*$-superrigid actions were only discovered very recently in my joint paper with Popa \cite{PV09}. We found a family of amalgamated free product groups $\Gamma$ and a large class of W$^*$-superrigid $\Gamma$-actions, including (generalized) Bernoulli actions, Gaussian actions and certain co-induced actions.

Note that W$^*$-superrigidity for an action $\Gamma \actson (X,\mu)$ is equivalent to the \lq sum\rq\ between its OE superrigidity and the uniqueness, up to unitary conjugacy, of $\rL^\infty(X)$ as a group measure space Cartan subalgebra in $\rL^\infty(X) \rtimes \Gamma$. This makes W$^*$-superrigidity results extremely difficult to obtain, since each one of these problems is notoriously hard. Contrary to the long list of OE superrigid actions referred to in the previous paragraph, unique Cartan decomposition proved to be much more challenging to establish, and the only existing results cover very particular group actions. Thus, a first such result, obtained by Ozawa and Popa \cite{OP07}, shows that given any profinite action $\Gamma \actson X$ of a product of free groups $\Gamma = \F_{n_1}\times \cdots \times \F_{n_k}$, with $k\geq 1$, $2\leq n_i \leq \infty$, all Cartan subalgebras of $M=\rL^\infty(X) \rtimes \Gamma$ are unitarily conjugate to $\rL^\infty(X)$. A similar result, covering groups $\Gamma$ that have the complete metric approximation property and that admit a proper $1$-cocycle into a non-amenable representation, was then proved in \cite{OP08}. More recently, Peterson showed \cite{Pe09} that factors arising from profinite actions of non-trivial free products
$\Gamma=\Gamma_1 * \Gamma_2$, with at least one of the $\Gamma_i$ not having the Haagerup property, have a unique group measure space Cartan
subalgebra, up to unitary conjugacy. But so far, none of these group actions could be shown to be OE superrigid. Nevertheless, an intricate combination of results in \cite{Io08,OP08,Pe09} were used to prove the existence of virtually\footnote{Following \cite{Fu98a}, {\it virtual} means that the ensuing conjugacy of $\Gamma \actson X$ and $\Lambda \actson Y$ is up to finite index subgroups of $\Gamma, \Lambda$.} W$^*$-superrigid group actions $\Gamma \actson X$ in \cite{Pe09}, by a Baire category argument.

In \cite{PV09} we established a very general unique Cartan decomposition result, which allowed us to obtain a wide range of W$^*$-superrigid group actions. Thus, we first proved the uniqueness, up to unitary conjugacy, of the group measure space Cartan subalgebra in the II$_1$ factor given by an arbitrary free ergodic p.m.p.\ action of any group $\Gamma$ belonging to a large family $\cG$ of amalgamated free product groups. By combining this with Kida's OE superrigidity in \cite{Ki09}, we deduced that if $T_n < \PSL(n,\Z)$ denotes the group of upper triangular matrices in $\PSL(n,\Z)$, then any free mixing p.m.p.\ action of $\Gamma=\PSL(n,\Z)*_{T_n}
\PSL(n,\Z)$ is W$^*$-superrigid. In combination with \cite{Po05,Po06a}, we proved that for many groups $\Gamma$ in the family $\cG$, the Bernoulli actions of $\Gamma$ are W$^*$-superrigid.

Very recently, Ioana \cite{Io10} obtained the beautiful result that all Bernoulli actions of property (T) groups are W$^*$-superrigid.

Recall from paragraph \ref{par.defrig} Popa's strong rigidity theorem for Bernoulli actions and note the asymmetry in the formulation: there is a (rigidity) condition on the group $\Gamma$ and a (deformation) condition on the action $\Lambda \actson (Y,\eta)$. One of the novelties of \cite{PV09} is a transfer of rigidity principle showing that under W$^*$-equivalence of $\Gamma \actson X$ and $\Lambda \actson Y$, some of the rigidity properties of $\Gamma$ persist in the arbitrary unknown group $\Lambda$. Note however that property (T) itself is not stable under W$^*$-equivalence: there exist W$^*$-equivalent group actions such that $\Gamma$ has property (T) while $\Lambda$ has not (see Section \ref{sec.cartan}).

\smalltitle{Indecomposability results}

Related to the uni\-que\-ness problem of Cartan subalgebras obviously is the existence question. Voi\-cu\-lescu \cite{Vo95} proved that the free group factors admits no Cartan subalgebra, because the presence of a Cartan subalgebra forces the free entropy dimension of any generating set to be smaller or equal than $1$. Another application of Voiculescu's free entropy theory was given by Ge \cite{Ge96} who showed that the free group factors are prime: they cannot be written as the tensor product of two II$_1$ factors.

Using delicate C$^*$-algebra techniques, Ozawa \cite{Oz03} proved that for all icc word hyperbolic groups $\Gamma$ -- in particular when $\Gamma$ equals the free group $\F_n$ -- the group factor $\rL \Gamma$ is solid: the relative commutant $A' \cap \rL \Gamma$ of an arbitrary diffuse\footnote{A von Neumann algebra is called diffuse if it admits no minimal projections.} subalgebra is injective. Obviously non-hyperfinite solid II$_1$ factors, as well as all their non-hyperfinite subfactors, are prime. A combination of techniques from \cite{Po01,Oz03} then allowed Ozawa and Popa \cite{OP03} to introduce a family of II$_1$ factors that have an essentially unique tensor product decomposition into prime factors. Peterson's $\rL^2$-rigidity \cite{Pe06} -- a II$_1$ factor analogue for the vanishing of the first $\ell^2$-Betti number of a group -- as well as Bass-Serre rigidity for free products of von Neumann algebras \cite{IPP05,CH08} provided further examples of prime factors.

The free group factors have no Cartan subalgebra and are solid. In \cite{OP07} both properties are brought together and $\rL \F_n$ is shown to be strongly solid: the normalizer of an arbitrary diffuse abelian subalgebra is hyperfinite. Other examples of strongly solid II$_1$ factors were given in \cite{Ho09,HS09}.

\smalltitlezonder{Organization of the paper}

In so far as the above gave an overview of some recent developments, in the rest of the paper I present the main ideas behind a number of chosen topics: Popa's deformation/rigidity theory in Section \ref{sec.defrig}, computations of fundamental groups in Section \ref{sec.fundamental}, the (non-)uniqueness of Cartan subalgebras in Section \ref{sec.cartan} and W$^*$-superrigidity in Section~\ref{sec.superrigid}.

\section{Preliminaries}\label{sec.prelim}

\smalltitlezonder{Traces, II$_1$ factors and the Hilbert bimodule $\rL^2(M)$}

A \emph{finite trace} on a von Neumann algebra $M$ is a linear map $\tau : M \recht \C$ satisfying $\tau(xy) = \tau(yx)$ for all $x,y \in M$. We say that $\tau$ is \emph{positive} if $\tau(x) \geq 0$ for all positive operators $x \in M$. A positive trace $\tau$ is called \emph{faithful} if the equality $\tau(x^* x) = 0$ implies that $x = 0$. A positive trace $\tau$ is called \emph{normal} if $\tau$ is weakly continuous on the unit ball of $M$.

A \emph{II$_1$ factor} is a von Neumann algebra with trivial center that admits a non-zero finite positive trace $\tau$ and that is non-isomorphic to a matrix algebra $\M_n(\C)$. Normalizing $\tau$ such that $\tau(1) = 1$, the trace is unique. Moreover $\tau$ is automatically normal. We denote by $\|x\|_2 = \sqrt{\tau(x^*x)}$ the $\rL^2$-norm corresponding to $\tau$. Completing $M$ w.r.t.\ the scalar product $\langle x,y\rangle = \tau(x^* y)$ yields the Hilbert space $\rL^2(M)$, which is an $M$-$M$-bimodule by left and right multiplication on $M$.

We denote by $\Tr$ the non-normalized trace on $\M_n(\C)$. Occasionally, $\Tr$ denotes the infinite trace on positive operators in $\B(\cH)$.

\smalltitlezonder{Group von Neumann algebras}

Let $\Gamma$ be a countable group. Then $\rL \Gamma$ is the unique tracial von Neumann algebra generated by unitary elements $(u_g)_{g \in \Gamma}$ with the following two properties: $u_g u_h = u_{gh}$ for all $g,h \in \Gamma$ and $\tau(u_g) = 0$ for all $g \neq e$. Alternatively, we denote by $(\delta_g)_{g \in \Gamma}$ the standard orthonormal basis of $\ell^2(\Gamma)$, define the translation unitary operators $u_g$ as $u_g \delta_h = \delta_{gh}$ and define $\rL \Gamma$ as the von Neumann algebra generated by $\{u_g \mid g \in \Gamma\}$, with $\tau$ being given by $\tau(x) = \langle \delta_e , x \delta_e\rangle$ for all $x \in \rL \Gamma$.

\smalltitlezonder{Group measure space construction}

If $(P,\tau)$ is a tracial von Neumann algebra and $\Gamma \overset{\alpha}{\actson} P$ is an action of a countable group $\Gamma$ by trace preserving automorphisms $\al_g \in \Aut P$, the \emph{crossed product} $P \rtimes \Gamma$ is the unique tracial von Neumann algebra $(M,\tau)$ generated by a trace-preserving copy of $P$ and unitary elements $(u_g)_{g \in \Gamma}$ satisfying the following properties:
\begin{align*}
& u_g a u_g^* = \al_g(a)\;\;\text{for all}\;\; g \in \Gamma , a \in P \; , \quad u_g u_h = u_{gh}\;\;\text{for all}\;\; g,h \in \Gamma \; , \\ &\tau(au_g) = 0 \;\;\text{for all}\;\; a \in P, g \neq e \; .
\end{align*}
The map $au_g \mapsto a \ot \delta_g$ provides an identification $\rL^2(P \rtimes \Gamma) = \rL^2(P) \ovt \ell^2(\Gamma)$ and then an explicit realization of $P \rtimes \Gamma$ as an algebra of bounded operators on the Hilbert space $\rL^2(P) \ovt \ell^2(\Gamma)$.

When $\Gamma \actson (X,\mu)$ is a probability measure preserving (p.m.p.) action, one considers the corresponding trace preserving action $\Gamma \actson \rL^\infty(X)$ and constructs $M = \rL^\infty(X) \rtimes \Gamma$. The abelian subalgebra $\rL^\infty(X) \subset M$ is maximal abelian if and only if $\Gamma \actson (X,\mu)$ is essentially free, meaning that for all $g \neq e$ the set $\{x \in X \mid g \cdot x = x\}$ has measure zero. When $\Gamma \actson (X,\mu)$ is essentially free, the center of $M$ equals $\rL^\infty(X)^\Gamma$, the algebra of $\Gamma$-invariant functions in $\rL^\infty(X)$. Hence, factoriality of $M$ is then equivalent with ergodicity of $\Gamma \actson (X,\mu)$.

\smalltitlezonder{(Generalized) Bernoulli actions}

If $\Gamma$ is an infinite countable group and if $(X_0,\mu_0)$ is a non-trivial probability space, define the infinite product $(X,\mu) := (X_0,\mu_0)^\Gamma$ on which $\Gamma$ acts by shifting the indices: $(g \cdot x)_h = x_{g^{-1}h}$. The action $\Gamma \actson (X,\mu)$ is called the \emph{Bernoulli action} with base space $(X_0,\mu_0)$ and it is a free ergodic p.m.p.\ action.

More generally, if $\Gamma$ acts on the countably infinite set $I$, one considers the \emph{generalized Bernoulli action} $\Gamma \actson (X_0,\mu_0)^I$ given by $(g \cdot x)_i = x_{g^{-1} \cdot i}$. This action is p.m.p.\ and it is ergodic if and only if every orbit $\Gamma \cdot i$ is infinite. If $(X_0,\mu_0)$ is non-atomic, essential freeness is equivalent with every $g \neq e$ acting non-trivially on $I$. If $(X_0,\mu_0)$ has atoms, essential freeness is equivalent with every $g \neq e$ moving infinitely many $i \in I$.

\smalltitlezonder{Bimodules}

An $M$-$N$-bimodule $\bim{M}{\cH}{N}$ between von Neumann algebras $M$ and $N$ is a Hilbert space $\cH$ equipped with a normal unital $*$-homomorphism $\lambda : M \recht \B(\cH)$ and a normal unital $*$-anti-homomorphism $\rho : N \recht \B(\cH)$ such that $\lambda(M)$ and $\rho(N)$ commute. We write $x  \xi y$ instead of $\lambda(x)\rho(y)\xi$.

Bimodules should be considered as the II$_1$ factor analogue of unitary group representations. Based on this philosophy, several representation theoretic properties of groups have a II$_1$ factor counterpart: amenability, the Haagerup property, property (T), etc.

To establish the dictionary between group representations and bimodules, take a countable group $\Gamma$ and put $M = \rL \Gamma$. Whenever $\pi : \Gamma \recht \cU(\cK)$ is a unitary representation, define the Hilbert space $\cH_\pi = \ell^2(\Gamma) \ovt \cK$ and turn $\cH_\pi$ into an $M$-$M$-bimodule by putting $u_g (\delta_h \ot \xi) u_k := \delta_{ghk} \ot \pi(g) \xi$. The trivial representation corresponds to the \emph{trivial bimodule} $\bim{M}{\rL^2(M)}{M}$, the regular representation corresponds to the \emph{coarse bimodule} $\bim{M \ot 1}{\bigl(\rL^2(M) \ovt \rL^2(M))}{1 \ot M}$ and one defines the notions of \emph{containment} and \emph{weak containment} of bimodules \cite{Po86} in such a way that through the construction $\pi \rightsquigarrow \cH_\pi$ these notions exactly correspond to the well known concepts from representation theory. Finally, the \emph{Connes tensor product} of bimodules \cite[V.Appendix B]{Co94} is so that $\cH_\pi \ot_M \cH_\rho \cong \cH_{\pi \otimes \rho}$.

\begin{definition}\label{def.amenT}
A tracial von Neumann algebra $(M,\tau)$ is called \emph{amenable} \cite{Po86} if the coarse bimodule weakly contains the trivial bimodule.

We say that $(M,\tau)$ has \emph{property (T)} \cite{CJ83} if any $M$-$M$-bimodule weakly containing the trivial bimodule, must contain the trivial bimodule.

We finally say that the subalgebra $N \subset M$ has the \emph{relative property (T)} \cite{Po01} if any $M$-$M$-bimodule weakly containing the trivial bimodule, must contain the bimodule $\bim{N}{\rL^2(M)}{M}$.
\end{definition}

\smalltitlezonder{Completely positive maps and bimodules}

A linear map $\vphi : M \recht N$ is called \emph{completely positive} if for every $n$ the amplified map $\id \ot \vphi : \M_n(\C) \ot M \recht \M_n(\C) \ot N$ maps positive operators to positive operators. In the same way as unitary representations are related to positive definite functions, also bimodules and completely positive maps form two sides of the same story. Whenever $\vphi : M \recht N$ is a unital trace preserving completely positive map, the separation and completion of the algebraic tensor product $M \otalg N$ w.r.t.\ the scalar product $\langle a \ot b, c \ot d\rangle = \tau(b^* \vphi(a^* c) d)$ defines a Hilbert space $\cH$ that naturally becomes an $M$-$N$-bimodule. By construction $\xi = 1 \ot 1$ is a cyclic vector, meaning that $M \xi N$ is dense in $\cH$, and is a trace vector, meaning that $\tau(a) = \langle \xi,a\xi\rangle$ and $\tau(b) = \langle \xi , \xi b \rangle$ for all $a \in M, b \in N$. Every bimodule with a cyclic trace vector arises in this way.

\smalltitlezonder{Cartan subalgebras and equivalence relations}

Whenever $A \subset M$ is a von Neumann subalgebra, we denote by $\cN_M(A) := \{u \in \cU(M) \mid uAu^* = A\}$ the group of unitaries in $M$ that normalize $A$. A \emph{Cartan subalgebra} $A \subset M$ of a II$_1$ factor is a maximal abelian subalgebra such that $\cN_M(A)$ generates $M$. Whenever $\Gamma \actson (X,\mu)$ is a free ergodic p.m.p.\ action, $\rL^\infty(X) \subset \rL^\infty(X) \rtimes \Gamma$ is a Cartan subalgebra. We call such Cartan subalgebras of \emph{group measure space} type.

The relevance of Cartan subalgebras in the study of group measure space factors stems from the following theorem.

\begin{theorem}[Singer \cite{Si55}]\label{thm.singer}
Let $\Gamma \actson (X,\mu)$ and $\Lambda \actson (Y,\eta)$ be free ergodic p.m.p.\ actions and denote $A = \rL^\infty(X)$, $B = \rL^\infty(Y)$. Assume that $\Delta : X \recht Y$ is an isomorphism of probability spaces with corresponding trace preserving isomorphism $\Delta_* : A \recht B : F \mapsto F \circ \Delta^{-1}$. Then, the following are equivalent.
\begin{itemize}
\item $\Delta$ is an orbit equivalence: for almost every $x \in X$, we have $\Delta(\Gamma \cdot x) = \Lambda \cdot \Delta(x)$.
\item $\Delta_*$ extends to a $*$-isomorphism $A \rtimes \Gamma \recht B \rtimes \Lambda$.
\end{itemize}
\end{theorem}

A \emph{II$_1$ equivalence relation} \cite{FM75} on a standard probability space $(X,\mu)$ is an equivalence relation $\cR \subset X \times X$ with \emph{countable} equivalence classes such that $\cR$ is a Borel subset of $X \times X$ and such that $\cR$ is \emph{ergodic} and \emph{probability measure preserving.} Here $\cR$ is called ergodic if every $\cR$-saturated Borel set has measure $0$ or $1$, while $\cR$ is said to be p.m.p.\ if every bimeasurable bijection $\vphi: X \recht X$ with graph inside $\cR$, preserves $\mu$.

Whenever $\Gamma \actson (X,\mu)$ is an ergodic p.m.p.\ action, the \emph{orbit equivalence relation} $\cR(\Gamma \actson X)$ is of type II$_1$. Every II$_1$ equivalence relation is of this form, but it is not always possible to choose an essentially \emph{free} action implementing $\cR$ (see \cite{Fu98b} and \cite[Section 7]{PV08b}).

A variant of the group measure space construction \cite{FM75} allows to associate a II$_1$ factor $\rL \cR$ to any II$_1$ equivalence relation $\cR$ on $(X,\mu)$. By construction $\rL \cR$ contains a copy of $\rL^\infty(X)$ as a Cartan subalgebra. Every Cartan inclusion $A \subset M$ arises in this way, modulo the possible appearance of a scalar $2$-cocycle on $\cR$. When $\Gamma \actson (X,\mu)$ is a free ergodic p.m.p.\ action, we canonically have $\rL^\infty(X) \rtimes \Gamma = \rL(\cR(\Gamma \actson X))$. It is however important to note that both the crossed product and the orbit equivalence relation make sense for non-free actions, but no longer yield isomorphic von Neumann algebras.

\smalltitlezonder{Fundamental groups of II$_1$ factors}

When $M$ is a II$_1$ factor and $t > 0$, one defines as follows the amplification $M^t$. For $0 < t \leq 1$, take a projection $p \in M$ with $\tau(p) = t$ and put $M^t := p M p$. For larger $t$, take an integer $n$, a projection $p \in \M_n(\C) \ot M$ with $(\Tr \ot \tau)(p) = t$ and put $M^t := p(\M_n(\C) \ot M)p$. As such, $M^t$ is well defined up to isomorphism. One proves that $(M^t)^s \cong M^{ts}$. The \emph{fundamental group} $\cF(M)$ is defined as the set of $t > 0$ such that $M^t \cong M$. The fundamental group of a II$_1$ equivalence relation $\cR$ is defined in a similar way. If $M = \rL \cR$ denotes the associated II$_1$ factor, by construction $\cF(\cR) \subset \cF(M)$, but this inclusion can be strict \cite[Section 6.1]{Po06a}.

\section{Popa's deformation/rigidity theory}\label{sec.defrig}

Popa's deformation/rigidity theory, initiated in \cite{Po01}, has revolutionized our understanding of II$_1$ factors. We explain in this section what kind of deformations Popa introduced and how they can be combined with the rigidity given by the relative property (T).

\begin{definition}\label{def.deform}
A \emph{deformation of the identity} on a tracial von Neumann algebra $(M,\tau)$ is a sequence of normal completely positive maps $\vphi_n : M \recht M$ that are unital, trace preserving and satisfy
$$\|\vphi_n(x) - x \|_2 \recht 0 \quad\text{for all}\;\; x \in M \; .$$
\end{definition}

Both group factors $\rL \Gamma$ and crossed products $P \rtimes \Gamma$ admit natural deformations of the identity.
Indeed, if $\vphi : \Gamma \recht \C$ is a positive definite function, both
\begin{align*}
& \rL \Gamma \recht \rL \Gamma : u_g \mapsto \vphi(g) u_g \quad\text{for all}\;\; g \in \Gamma \quad\text{and}\\ & P \rtimes \Gamma \recht P \rtimes \Gamma : a u_g \mapsto \vphi(g) a u_g \quad\text{for all}\;\; a \in P, g \in \Gamma
\end{align*}
extend to normal completely positive maps on $\rL \Gamma$, resp.\ $P \rtimes \Gamma$. If $\vphi$ is normalized, i.e.\ $\vphi(e) = 1$, these maps are unital and trace preserving.

\begin{example} \label{ex.Haag-freeprod}
If $\Gamma$ has the Haagerup approximation property \cite{Ha78}, there exists a sequence $\vphi_n : \Gamma \recht \C$ of positive definite functions converging pointwise to $1$ and with $\vphi_n \in c_0(\Gamma)$ for every $n$. As we discuss below, the corresponding deformation of the identity of $P \rtimes \Gamma$ plays a crucial role in \cite{Po01}.

If $\Gamma = \Gamma_1 * \Gamma_2$ is a free product and $|g|$ denotes the natural word length of $g \in \Gamma$ w.r.t.\ this free product decomposition and if $0 < \rho < 1$, then the formula $\vphi_\rho(g) = \rho^{|g|}$ defines a positive definite function on $\Gamma$. If $\rho \recht 1$, then $\vphi_\rho \recht 1$ pointwise. The corresponding deformation of the identity is the starting point for \cite{IPP05,PV09} and also this is discussed below.
\end{example}

A second family of deformations of the identity arises as follows. Let $\Gamma \overset{\al}{\actson} (P,\tau)$ be a trace preserving action. Assume that $\vphi_n : P \recht P$ is a deformation of the identity such that $\vphi_n \circ \al_g = \al_g \circ \vphi_n$ for all $g \in \Gamma$, $n \in \N$. Then, the formula $a u_g \mapsto \vphi_n(a) u_g$ defines a deformation of the identity on $P \rtimes \Gamma$.

\begin{definition}[\cite{Po03}]
A p.m.p.\ action $\Gamma \actson (X,\mu)$ is \emph{malleable} if there exists a continuous family $(\al_t)_{t \in [0,1]}$ of p.m.p.\ transformations of $X \times X$ such that for all $t \in [0,1]$, $\al_t$ commutes with the diagonal action $g \cdot (x,y) = (g \cdot x,g \cdot y)$ and such that $\al_0 = \id$ and $\al_1(x,y)$ is of the form $(\ldots,x)$.
\end{definition}

The \emph{Bernoulli action} $\Gamma \actson [0,1]^\Gamma$ is malleable. It suffices to construct a continuous family $(\al^0_t)_{t \in [0,1]}$ of p.m.p.\ transformations of the square $[0,1] \times [0,1]$ such that $\al^0_0 = \id$ and $\al^0_1(x,y) = (\ldots,x)$. This can be done by \lq rotating the square\rq\ counterclockwise over 90 degrees. Next, one identifies $[0,1]^\Gamma \times [0,1]^\Gamma = ([0,1] \times [0,1])^\Gamma$ and defines $(\al_t(x,y))_g = \al^0_t(x_g,y_g)$. By construction, $\al_t$ commutes with the diagonal $\Gamma$-action.

Other examples of malleable actions are \emph{generalized Bernoulli actions} $\Gamma \actson [0,1]^I$ given by an action $\Gamma \actson I$ or \emph{Gaussian actions} given by an orthogonal representation of $\Gamma$ on a real Hilbert space.

\begin{example} \label{ex.malleable}
If $\al_t$ is a malleable deformation, put $A = \rL^\infty(X)$ and define the corresponding automorphisms $\al_t$ of $A \ovt A = \rL^\infty(X \times X)$ given by $\al_t(F) := F(\al_t^{-1}(\, \cdot \,))$. By definition $\al_0 = \id$ and $\al_1(a \ot 1) = 1 \ot a$ for all $a \in A$. View $A \hookrightarrow A \ovt A : a \mapsto a \ot 1$ and denote by $E : A \ovt A \recht A$ the trace preserving conditional expectation (which corresponds to integration w.r.t.\ the second variable). The formula $\vphi_t : A \recht A : \vphi_t(a) = E_A(\al_t(a \ot 1))$ defines a continuous family of unital trace preserving completely positive maps with $\vphi_0 = \id$ and $\vphi_1(a) = \tau(a)1$. When $t \recht 0$, we get a deformation of the identity on $\rL^\infty(X) \rtimes \Gamma$ that is at the heart of \cite{Po03,Po04}.
\end{example}

A variant of the malleable deformation for Bernoulli actions is the \emph{tensor length deformation} \cite{Io06}. Indeed, given a base probability space $(X_0,\mu_0)$ and a countable set $I$, put $(X,\mu) = (X_0,\mu_0)^I$ and identify $A:=\rL^\infty(X)$ with the infinite tensor product $\overline{\otimes}_{i \in I} (A_0,\tau)$, where $A_0 := \rL^\infty(X_0)$. We write $A = A_0^I$. Whenever $J \subset I$ is a subset, we view $A_0^J$ as a subalgebra of $A_0^I$. We then define for every $0 < \rho < 1$,
$$\theta_\rho : A \recht A : \theta_\rho(a) = \rho^n a \quad\text{when}\;\; a \in (A_0 \ominus \C 1)^J \;\;\text{and}\;\; |J| = n \; .$$
Then $\theta_\rho$ is a well defined unital trace preserving normal completely positive map and $\theta_\rho \recht \id$ when $\rho \recht 1$. By construction, $\theta_\rho$ commutes with the generalized Bernoulli action $\Gamma \actson A_0^I$ whenever $\Gamma \actson I$.

\smalltitlezonder{Combining deformation and rigidity}

Let $N \subset M$ be an inclusion with the relative property (T), see Definition \ref{def.amenT}. Whenever $\vphi_n : M \recht M$ is a deformation of the identity, it follows that $\vphi_n$ converges uniformly in $\|\,\cdot\,\|_2$ on the unit ball of $N$.
We illustrate this combination of deformation and rigidity by indicating the main ideas behind two of Popa's theorems.

\begin{theorem}[Popa \cite{Po01}] \label{thm.popatrivfund}
The II$_1$ factor $M = \rL(\Z^2 \rtimes \SL(2,\Z))$ has trivial fundamental group.
\end{theorem}

Write $A = \rL^\infty(\T^2)$, $\Gamma = \SL(2,\Z)$ and identify $M = A \rtimes \Gamma$. The group $\Gamma$ has the Haagerup property and hence admits a sequence of positive definite functions $\vphi_n$ such that $\vphi_n \in c_0(\Gamma)$ for all $n$ and $\vphi_n \recht 1$ pointwise. As in Example \ref{ex.Haag-freeprod}, we get a deformation of the identity on $M$ given by $\theta_n(a u_g) = \vphi_n(g) a u_g$.
Assume that $N \subset M$ has the relative property (T) and choose $\eps > 0$. Because $\vphi_n \in c_0(\Gamma)$, we find $n \in \N$ and a finite subset $\cF \subset \Gamma$ such that for all $b$ in the unit ball of $N$, the $\|\,\cdot\,\|_2$-distance of $b$ to $\lspan \{a u_g \mid a\in A, g \in \cF \}$ is smaller than $\eps$. The only \lq obvious\rq\ subalgebras of $A \rtimes \Gamma$ with such an approximation property are those that are unitarily conjugate to a subalgebra of $A$, i.e.\ $v N v^* \subset A$ for some $v \in \cU(M)$. Popa's \emph{intertwining-by-bimodules,} that we recall below, ensures that this feeling is indeed (almost) correct. It is even exactly correct when moreover $N \subset M$ is a Cartan subalgebra.

Now observe that $A \subset M$ has the relative property (T) and is a Cartan subalgebra. So, whenever $\al$ is an automorphism of $M$, the subalgebra $\al(A)$ still has the relative property (T) and the previous paragraph implies that, up to a unitary conjugacy, every automorphism of $M$ globally preserves $A$. This means that every automorphism of $M$ induces an automorphism of the orbit equivalence relation of $\SL(2,\Z) \actson \T^2$. A similar statement is true for isomorphisms $M \recht pMp$ and therefore, the fundamental group of $M$ equals the fundamental group of the equivalence relation, which is trivial by \cite{Ga99,Ga01}.

\begin{theorem}[Popa \cite{Po03,Po04}] \label{thm.strongrigid}
Let $\Gamma$ be a property (T) group and $\Gamma \actson (X,\mu)$ any free ergodic p.m.p.\ action. Let $\Lambda$ be any icc group and $\Lambda \actson (Y_0,\eta_0)^\Lambda$ the Bernoulli action. Put $(Y,\eta) := (Y_0,\eta_0)^\Lambda$.

If $\rL^\infty(X) \rtimes \Gamma \cong \rL^\infty(Y) \rtimes \Lambda$, then the groups $\Gamma,\Lambda$ are isomorphic and their actions are conjugate.
\end{theorem}

Put $A = \rL^\infty(X)$ and $B_0 = \rL^\infty(Y_0)$. Assume that $A \rtimes \Gamma = B_0^\Lambda \rtimes \Lambda$ and consider on $B_0^\Lambda \rtimes \Lambda$ the tensor length deformation $\theta_\rho$ defined after Example \ref{ex.malleable}. Since the subalgebra $\rL \Gamma$ has property (T), for $\rho$ close enough to $1$, we get that $\theta_\rho$ is uniformly close to the identity on the unit ball of $\rL \Gamma$. When $b \in B_0^\Lambda$, the norm $\|\theta_\rho(b) - b \|_2$ is small when $b$ can be written as a linear combination of \lq short\rq\ elementary tensors. The only obvious \lq short\rq\ subalgebras of $B_0^\Lambda \rtimes \Lambda$ are those that can be unitarily conjugated into either $\rL \Lambda$ or $B_0^J$ for some finite subset $J \subset \Lambda$. The abelian algebra $B_0^J$ can never house the property (T) algebra $\rL \Gamma$ and Popa indeed manages to prove that $\rL \Gamma$ must be unitarily conjugate to a subalgebra of $\rL \Lambda$.

So we may assume that $\rL \Gamma \subset \rL \Lambda$. In a second and analytically very delicate part, Popa basically proves the following: if a subalgebra $A$ of $B_0^\Lambda \rtimes \Lambda$ is both abelian and normalized by many unitaries in $\rL \Lambda$, then $A$ must be unitarily conjugate to a subalgebra of $B_0^\Lambda$.

This brings us in the situation where $A$ can be unitarily conjugated into $B_0^\Lambda$ and $\rL \Gamma$ can be unitarily conjugated into $\rL \Lambda$. Popa proves that automatically both unitary conjugations can be done with the same unitary, yielding isomorphism of the groups and conjugacy of the actions.

\smalltitlezonder{Popa's intertwining-by-bimodules}

In \cite{Po01,Po03} Popa developed a powerful technique to approach the following question: when are two subalgebras $N,P \subset M$ unitarily conjugate? A detailed explanation and motivation for his method can be found in \cite[Section 5]{Po06b} and \cite[Appendix F]{BO08}. So, we are rather brief here. First consider the case where $M = P \rtimes \Gamma$ for some trace preserving action $\Gamma \actson P$. Every element $x \in M$ has a unique \emph{Fourier expansion}
$$x = \sum_{g \in \Gamma} x_g u_g \quad\text{with}\;\; x_g \in P$$
converging in $\|\,\cdot\,\|_2$. We call the $x_g \in P$ the Fourier coefficients of $x$.

\begin{theorem}[Popa \cite{Po01,Po03}]\label{thm.intertwine}
Let $N \subset P \rtimes \Gamma$ be a von Neumann subalgebra. Then the following two conditions are equivalent.
\begin{itemize}
\item There exist projections $p \in P, q \in N$, a normal unital $*$-homomorphism $\vphi : qNq \recht pPp$ and a non-zero partial isometry $v \in q (P \rtimes \Gamma) p$ satisfying $a v = v\vphi(a)$ for all $a \in qNq$.

\item There is no sequence of unitaries $v_n \in N$ whose Fourier coefficients converge to $0$ pointwise in $\|\,\cdot\,\|_2$, i.e.\ $\|(v_n)_g\|_2 \recht 0$ for all $g \in \Gamma$.
\end{itemize}
\end{theorem}

When $P \subset M$ is no longer the \lq core\rq\ of a crossed product, there is no notion of Fourier coefficients and their convergence to $0$ has to be replaced by the condition $\|E_P(x v_n y)\|_2 \recht 0$ for all $x,y \in M$, where $E_P : M \recht P$ is the unique trace preserving conditional expectation. If $M = P \rtimes \Gamma$, note that $(v_n)_g = E_P(v_n u_g^*)$.

The first condition in Theorem \ref{thm.intertwine} is of course not saying that $v$ is a unitary satisfying $v^* N v \subset P$. The left support projection of $v$ lies in the relative commutant of $qNq$ -- which in concrete applications is usually known -- but the right support projection of $v$ lies in the relative commutant of $\vphi(qNq)$ which is of course a priori unknown, since we do not know $\vphi$. Several techniques based on \emph{mixing properties} have been developed to take care of this relative commutant issue and we refer to \cite[Section 5]{Po06b} for a more detailed explanation.

\section{Fundamental groups of II$_1$ factors}\label{sec.fundamental}

In order to construct II$_1$ factors with a prescribed fundamental group $\cF$, you first need a cute construction of a II$_1$ factor $M$ such that all $t \in \cF$ obviously belong to $\cF(M)$ and then you need a powerful theory to make sure that no other $t > 0$ belong to $\cF(M)$. As an illustration, we first briefly explain two constructions that produce II$_1$ factors with prescribed countable fundamental group.

\smalltitlezonder{Connes-St{\o}rmer Bernoulli actions \cite{Po03}}

Let $(X_0,\mu_0)$ be an atomic probability space and $\Gamma$ a countable group. Put $(X,\mu) = (X_0,\mu_0)^\Gamma$ and define on $(X,\mu)$ the following II$_1$ equivalence relation: $x \sim y$ if and only if there exists a $g \in \Gamma$ and a finite subset $J \subset \Gamma$ such that $x_{gh} = y_h$ for all $h \in \Gamma - J$ and such that $\prod_{h \in J} \mu_0(x_{gh}) = \prod_{h \in J} \mu_0(y_h)$. Whenever $a \in X_0$, define the subset $Y_a \subset X$ as $Y_a := \{x \in X \mid x_e = a\}$. Given $a,b \in X_0$, the map $Y_a \recht Y_b$ changing $x_e$ from $a$ to $b$ and leaving the other $x_g$ untouched is an isomorphism of the restricted equivalence relations. Hence, $\mu_0(a) / \mu_0(b)$ belongs to $\cF(\cR)$ for all $a, b \in X_0$.

Denote by $M = \rL \cR$ the II$_1$ factor associated with $\cR$ and denote by $\cF$ the subgroup of $\R^*_+$ generated by all the ratios
$\mu_0(a) / \mu_0(b)$. We have seen that $\cF \subset \cF(M)$. In \cite{Po03} Popa shows that taking $\Gamma = \SL(2,\Z) \ltimes \Z^2$, one has the equality $\cF = \cF(M)$, with one of the ingredients of the proof being the triviality of the fundamental group of $\rL \Gamma$.

\smalltitlezonder{Free products of amplifications \cite{IPP05}}

Recall from Section \ref{sec.prelim} the notation $M^t$ for the amplification of a II$_1$ factor $M$.
In \cite{DR98} Dykema and R\u{a}dulescu established the following remarkable formula for an infinite free product of II$_1$ factors $M_n$~:
$$\Bigl( \underset{n \in \N}{*} M_n \Bigr)^t \cong \underset{n \in \N}{*} M_n^t \; .$$
So, whenever $\cF \subset \R^*_+$ is a countable subgroup different from $\{1\}$ and whenever $M$ is a II$_1$ factor, we put $P = *_{t \in \cF} M^t$ and conclude that both $\cF$ and $\cF(M)$ are subgroups of $\cF(P)$. As a consequence of their Bass-Serre theory for free products of II$_1$ factors, Ioana, Peterson and Popa \cite{IPP05} manage to prove that $\cF = \cF(P)$ when you take $M = \rL (\SL(2,\Z) \ltimes \Z^2)$. A more elementary variant of the previous construction was proposed in \cite{Ho07}.

Both constructions presented so far can be used to produce II$_1$ factors with arbitrary fundamental group $\cF$. However, for uncountable subgroups $\cF \subset \R^*_+$ the resulting II$_1$ factors do not have separable predual.

\smalltitlezonder{II$_1$ factors with uncountable fundamental group \cite{PV08a,PV08c}}

Again we start with a construction of a II$_1$ factor sowing a number of elements into the fundamental group. Let $\Gamma \actson (Z,\gamma)$ be a free ergodic action preserving the \emph{infinite} non-atomic measure $\gamma$. Put $N := \rL^\infty(Z) \rtimes \Gamma$ and note that $N$ is a II$_\infty$ factor. Whenever $\Delta$ is a non-singular\footnote{This means that $\Delta$ preserves sets of measure zero.} automorphism of $(Z,\gamma)$ satisfying $\Delta(g \cdot z) = g \cdot \Delta(z)$ a.e., the Radon-Nikodym derivative between $\gamma \circ \Delta^{-1}$ and $\gamma$ is $\Gamma$-invariant and hence constant a.e. So, $\Delta$ scales the infinite measure $\gamma$ by a factor that we denote by $\module \Delta$. Denote the group of all these automorphisms $\Delta$ as $\Centr_{\Aut Z}(\Gamma)$. Every $\Delta \in \Centr_{\Aut Z}(\Gamma)$ gives rise to an automorphism $\Delta_*$ of $N$ scaling the trace with the same factor $\module \Delta$. Whenever $p \in N$ is a projection of finite trace, put $M = pNp$ and note that $M$ is a II$_1$ factor. By construction we have
\begin{equation}\label{eq.inclu}
\module \Centr_{\Aut Z}(\Gamma) \subset \cF(pNp) \; .
\end{equation}
In the following paragraphs we explain the three main aspects of \cite{PV08a,PV08c}: how to get an equality in \eqref{eq.inclu}, how to make sure that $pNp$ is itself a group measure space II$_1$ factor and finally, how wild $\module \Centr_{\Aut Z}(\Gamma)$ can be. At the end this will give a feeling for the validity of the following theorem.

\begin{theorem}[Popa, Vaes \cite{PV08a,PV08c}]
Let $\Gamma_0$ be a non-trivial group and $\Sigma$ an infinite amenable group. Put $\Gamma = \Gamma_0^{*\infty}*\Sigma$. There exists a free ergodic p.m.p.\ action $\Gamma \actson (X,\mu)$ such that $\rL^\infty(X) \rtimes \Gamma$ has a fundamental group of arbitrary prescribed Hausdorff dimension between $0$ and $1$.

Actually, $\cF(\rL^\infty(X) \rtimes \Gamma)$ can be any group in the family
\begin{align*}
\Scentr := \{ \cF \subset \R^*_+ \mid \; & \text{there exists an amenable $\Lambda$ and a free ergodic}\;\; \Lambda \actson (Y,\eta) \\ &\text{preserving $\eta$ such that}\;\; \cF = \module \Centr_{\Aut Y}(\Lambda) \; \}
\end{align*}
\end{theorem}

\smalltitlezonder{How to get equality in \eqref{eq.inclu}}

Take $\Gamma$ of the form $\Gamma_1 * \Lambda$ and assume that $\Gamma \actson (X,\mu)$ is a free p.m.p.\ action with $\Gamma_1$ acting ergodically. Let $\Lambda \actson (Y,\eta)$ be a free ergodic action preserving the infinite non-atomic measure $\eta$. Put $(Z,\gamma) = (X \times Y,\mu \times \eta)$ and consider the action $\Gamma \actson Z$ given by $g \cdot (x,y) = (g \cdot x,y)$ if $g \in \Gamma_1$ and $h \cdot (x,y) = (h \cdot x, h \cdot y)$ if $h \in \Lambda$. It is easy to see that $\module \Centr_{\Aut Z}(\Gamma)$ equals $\module \Centr_{\Aut Y}(\Lambda)$. We now make the following assumptions on the action $\Gamma_1 * \Lambda \actson (X,\mu)$.
\begin{enumerate}
\item The action $\Gamma_1 \actson (X,\mu)$ is rigid, i.e.\ $\rL^\infty(X) \subset \rL^\infty(X) \rtimes \Gamma_1$ has the relative property (T) in the sense of Definition \ref{def.amenT}.
\item The group $\Lambda$ is amenable.
\item If $\phi : X_0 \recht X_1$ is a non-singular isomorphism between the non-negligible subsets $X_0,X_1 \subset X$ satisfying $\phi(X_0 \cap \Gamma_1 \cdot x) \subset \Gamma \cdot \vphi(x)$ for a.e.\ $x \in X_0$, then $\phi(x) \in \Gamma \cdot x$ for a.e.\ $x \in X_0$.
\end{enumerate}
For the following heuristic reasons, these assumptions imply that the inclusion in \eqref{eq.inclu} actually is an equality. Consider the II$_\infty$ factor $N = \rL^\infty(X \times Y) \rtimes (\Gamma_1*\Lambda)$ and let $\al$ be an automorphism of $N$. The Bass-Serre theory of \cite{IPP05} and the relative property (T) of $\rL^\infty(X)$ imply that after a unitary conjugacy, $\al$ globally preserves $\rL^\infty(X \times Y)$. In a next step, the relative property (T) of $\rL^\infty(X)$ together with the amenability of $\Lambda$ roughly implies that $\al$ globally preserves $\rL^\infty(X) \ot 1$. Then the third assumption above implies that we may assume that $\al(a \ot 1) = a \ot 1$ for all $a \in \rL^\infty(X)$ and hence that $\al$ is induced by an element of $\Centr_{\Aut Y}(\Lambda)$.

It is highly non-trivial to find group actions $\Gamma_1 * \Lambda \actson (X,\mu)$ satisfying the conditions 1, 2 and 3 above. Actually, in \cite{PV08a,PV08c} we use a Baire category argument to prove their existence whenever $\Gamma_1$ is an infinite free product $\Gamma_1 = \Gamma_0^{*\infty}$ with $\Gamma_0$ being an arbitrary non-trivial group. Replacing $\Gamma_0$ by $\Gamma_0 * \Gamma_0$, we may assume that $\Gamma_0$ is infinite. The proof roughly goes as follows: start with arbitrary free ergodic p.m.p.\ actions $\Gamma_0 \actson (X,\mu)$ and $\Lambda \actson (X,\mu)$ and view $\Gamma, \Lambda$ as subgroups of $\Aut(X,\mu)$. By \cite{Ga08}, there exists an automorphism $\beta_1 \in \Aut(X,\mu)$ such that the subgroups $\Gamma_0$ and $\beta_1 \Gamma_0 \beta_1^{-1}$ of $\Aut(X,\mu)$ generate a \emph{free and rigid} action of $\Gamma_0^{*2} := \Gamma_0 * \Gamma_0$ on $(X,\mu)$. By \cite{IPP05,To05} there exists $\psi \in \Aut(X,\mu)$ such that together with $\psi \Lambda \psi^{-1}$, we obtain a free action of $\Gamma_0^{*2} * \Lambda$ on $(X,\mu)$. We now start adding copies of $\Gamma_0$ acting as $\beta_n \Gamma_0 \beta_n^{-1} \subset \Aut(X,\mu)$ for well chosen $\beta_n \in \Aut(X,\mu)$. At stage $n$, given the free action $\Gamma_0^{*n} * \Lambda \actson X$, the rigidity of the action $\Gamma_0^{*2} \actson (X,\mu)$ implies that there are essentially only countably many partial isomorphisms $\phi$ that map $\Gamma_0^{*2}$-orbits into $\Gamma_0^{*n}*\Lambda$-orbits. By \cite{IPP05,To05} there exists $\beta_{n+1} \in \Aut(X,\mu)$ such that $\beta_{n+1}\Gamma_0 \beta_{n+1}^{-1}$ is free w.r.t.\ these countably many partial isomorphisms. Adding this new $\Gamma_0$-action, we obtain a free action of $\Gamma_0^{*(n+1)} * \Lambda$. Continuing by induction, we get the free action $\Gamma_0^{*\infty} * \Lambda \actson (X,\mu)$. We finally prove that there exists an infinite subset $E \subset \N$ containing $\{0,1\}$ such that $(*_{n \in E} \Gamma_0)*\Lambda \actson (X,\mu)$ satisfies condition 3. Since $\{0,1\} \subset E$ and $\Gamma_0^{*2} \actson (X,\mu)$ is rigid, condition 1 is satisfied as well.

\smalltitlezonder{Is $pNp$ itself a group measure space factor}

Take $\Gamma_1 * \Lambda \actson X \times Y$ as above, with $\Lambda$ being infinite amenable. Let $\Sigma$ be any infinite amenable group. By \cite[Lemma 3.6]{PV08c}, the restriction of the orbit equivalence relation $\cR(\Gamma_1 * \Lambda \actson X \times Y)$ to a set of finite measure is implemented by a free action of $\Gamma_1^{*\infty} * \Sigma$. So, $pNp$ is a group measure space factor.

\smalltitlezonder{How wild can $\module \Centr_{\Aut Y}(\Lambda)$ be}

In \cite{MNP68} the notion of an \emph{ergodic measure} $\nu$ on the real line is introduced: it is a $\sigma$-finite measure on the Borel sets of $\R$ such that for every $x \in \R$, the translation $\nu_x$ of $\nu$ by $x$ is either singular w.r.t.\ $\nu$ or equal to $\nu$, and such that denoting
$$H_\nu := \{x \in \R \mid \nu_x = \nu \}$$
the action of $H_\nu$ on $(\R,\nu)$ by translation is ergodic. More precisely, if $F : \R \recht \R$ is a Borel function and if for every $x \in H_\nu$ we have $F(x+y) = F(y)$ for $\nu$-a.e.\ $y$, then $F$ is constant $\nu$-a.e. The easiest examples of ergodic measures are the Lebesgue measure and the counting measure on a countable subgroup of $\R$.

In \cite{Aa86} it is shown that all subgroups of $\R^*_+$ of the form $\exp(H_\nu)$ arise as $\module \Centr_{\Aut Y}(\Z)$. Allowing more general amenable groups $\Lambda$ instead of $\Z$, this can be easily seen as follows. Assume that $H_\nu \neq \{0\}$. Viewing $H_\nu$ as a closed subgroup of $\Aut(\rL^\infty(\R,\nu))$ one turns $H_\nu$ into a Polish group. Take a countable dense subgroup $Q \subset H_\nu$ and define the additive subgroup $R \subset \R$ as $R = \Z[\exp(Q)]$. Then $Q$ acts on $R$ through multiplication by $\exp(q)$ and we put $\Lambda := R \rtimes Q$. Define the measure $\nutil$ on $\R$ such that $d\nutil(x) = \exp(-x) d\nu(x)$ and denote by $\lambda$ the Lebesgue measure on $\R$. Put $(Y,\eta) := (\R^2,\lambda \times \nutil)$. The action $\Lambda \actson Y$ given by $(r,q) \cdot (x,y) = (r + \exp(q)x, q+y)$ is measure preserving and ergodic. One checks easily that $\module \Centr_{\Aut Y}(\Lambda) = \exp(H_\nu)$.

It is also shown in \cite{Aa86} (cf.\ \cite[page 389]{PV08a}) that $H_\nu$ can be uncountable without being $\R^*_+$ and that actually $H_\nu$ can have any Hausdorff dimension between $0$ and $1$.

\begin{problem}
Give an intrinsic description of the subgroups of $\R^*_+$ that are of the form $\cF(M)$ where $M$ is a II$_1$ factor with separable predual.
\end{problem}

The only known a priori restriction on $\cF(M)$ is given by \cite[Proposition 2.1]{PV08a}: $\cF(M)$ is a Borel subset of $\R^*_+$ and carries a unique Polish group topology whose Borel sets are precisely the ones inherited from $\R^*_+$. It is however hard to believe that all \lq Polishable\rq\ Borel subgroups of $\R^*_+$ can arise as the fundamental group of a II$_1$ factor with separable predual.

\smalltitlezonder{Property (T) and fundamental groups}

Roughly speaking, the presence of property (T) forces fundamental groups to be countable. When $\Gamma$ is an icc property (T) group, Connes \cite{Co80} proved that $\cF(\rL \Gamma)$ is countable and the same method \cite{GG88} yields the countability of $\cF(\rL^\infty(X) \rtimes \Gamma)$ and $\cF(\cR(\Gamma \actson X))$, for all free ergodic p.m.p.\ actions $\Gamma \actson (X,\mu)$. For non-icc property (T) groups, Ioana \cite{Io08} could still prove the countability of $\cF(\cR(\Gamma \actson X))$, but in \cite{PV08c} we proved that if $\Gamma$ is a property (T) group whose virtual center\footnote{The virtual center is the set of elements with finite conjugacy class. It is a normal subgroup.} is not virtually abelian\footnote{A group is virtually abelian if it has an abelian subgroup of finite index.}, then $\Gamma$ admits a free ergodic p.m.p.\ action such that $\rL^\infty(X) \rtimes \Gamma$ is McDuff and in particular, has fundamental group $\R^*_+$. We were informed by Ershov that such groups exist as quotients of Golod-Shafarevich groups with property (T).

Zimmer \cite{Zi80} introduced a notion of property (T) for II$_1$ equivalence relations which is such that for free ergodic p.m.p.\ actions, the orbit equivalence relation $\cR(\Gamma \actson X)$ has property (T) if and only if the group $\Gamma$ has property (T). Using techniques from \cite{AD04} we proved \cite{PV08b} that for $n \geq 4$, the restriction of the II$_\infty$ relation $\cR(\SL(n,\Z) \actson \R^n)$ to a subset of finite Lebesgue measure, is a II$_1$ equivalence relation with property (T) in the sense of Zimmer and with fundamental group $\R^*_+$. In particular, this II$_1$ equivalence relation cannot be implemented by a free action of a group.

\section{(Non-)uniqueness of Cartan subalgebras}\label{sec.cartan}

\smalltitlezonder{Non-uniqueness of Cartan subalgebras}

Connes and Jones \cite{CJ81} have given the first examples of II$_1$ factors $M$ having more than one Cartan subalgebra up to conjugacy by an automorphism of $M$. Their construction goes as follows. Take a finite non-abelian group $\Sigma_0$, build the countable group $\Sigma = \Sigma_0^{(\N)}$ and the compact group $K = \Sigma_0^{\N}$ that we equip with its Haar measure. Consider the action $\Sigma \actson K$ by translation. Finally, let $\Gamma$ be any non-amenable group. Put $X:=K^\Gamma$ and consider the diagonal action $\Sigma \actson X$ which commutes with the Bernoulli action $\Gamma \actson X$. We obtain a free ergodic p.m.p.\ action $\Gamma \times \Sigma \actson X$ and hence, the Cartan subalgebra $A = \rL^\infty(X)$ of $M = \rL^\infty(X) \rtimes (\Gamma \times \Sigma)$. Taking non-commuting elements $g,h \in \Sigma_0$ and defining $g_n,h_n \in \Sigma$ as being $g,h$ in position $n$, one obtains two non-commuting central sequences in $M$. Hence, $M$ is a \emph{McDuff factor} \cite{McD70}, which means that $M \cong M \ovt R$, where $R$ denotes the hyperfinite II$_1$ factor. Choosing any Cartan subalgebra $B \subset R$, one transports back the Cartan subalgebra $A \ovt B \subset M \ovt R$ to a Cartan subalgebra of $M$ whose associated equivalence relation is not strongly ergodic. The initial equivalence relation given by $A \subset M$ is strongly ergodic, so that both Cartan subalgebras are non-conjugate by an automorphism of $M$.

Ozawa and Popa \cite{OP07} provided examples of II$_1$ factors where one can explicitly see two Cartan subalgebras. We first explain the general procedure and provide a concrete example below. Assume that $\Gamma$ is a countable group and $\Sigma \lhd \Gamma$ an infinite abelian normal subgroup. Assume that $\Sigma \lhd \Gamma$ has the following relative icc property: for every $g \in \Gamma - \Sigma$, the set $\{s g s^{-1} \mid s \in \Sigma \}$ is infinite. Let $\Sigma \hookrightarrow K$ be a dense embedding of $\Sigma$ into a compact abelian group $K$. Equip $K$ with its Haar measure and define the action $\Sigma \actson K$ by translation. Assume that we are given an extension of this action to an essentially free p.m.p.\ action $\Gamma \actson K$. Then,
$\rL^\infty(K)$ and $\rL \Sigma$ are non-unitarily conjugate Cartan subalgebras of $\rL^\infty(K) \rtimes \Gamma$.

An interesting concrete example is given by $\Gamma = \Z^n \rtimes \SL(n,\Z)$ and its natural affine action on $\Z_p^n$. We get
$$\rL^\infty(\Z_p^n) \rtimes (\Z^n \rtimes \SL(n,\Z)) = \rL^\infty(\T^n) \rtimes (\widehat{\Z_p^n} \rtimes \SL(n,\Z)) \; .$$
If $n = 2$, the group $\widehat{\Z_p^n} \rtimes \SL(n,\Z)$ has the Haagerup property while $\Z^n \rtimes \SL(n,\Z)$ does not. If $n \geq 3$, the group $\Z^n \rtimes \SL(n,\Z)$ has property (T) while $\widehat{\Z_p^n} \rtimes \SL(n,\Z)$ does not. So, neither property (T) nor the Haagerup property are stable under W$^*$-equivalence. Since they are stable under orbit equivalence, it follows that the Cartan subalgebras $\rL^\infty(\Z_p^n)$ and $\rL^\infty(\T^n)$ are non-conjugate by an automorphism of the ambient II$_1$ factor.

\smalltitlezonder{Uniqueness of Cartan subalgebras}

We mentioned in Section \ref{sec.panorama} that Ozawa and Popa \cite{OP07,OP08} established the first -- and up to now, only -- uniqueness results for Cartan subalgebras up to unitary conjugacy. Recall that a \emph{profinite} p.m.p.\ action $\Gamma \actson (X,\mu)$ is by definition the inverse limit $\invlimit (X_n,\mu_n)$ of a directed system of actions $\Gamma \actson (X_n,\mu_n)$ on \emph{finite} probability spaces.

\begin{theorem}[Ozawa, Popa \cite{OP07}] \label{thm.uniqueCartan}
Let $n \geq 2$ and let $\F_n \actson (X,\mu)$ be an ergodic profinite p.m.p.\ action. Put $A = \rL^\infty(X)$ and $M = A \rtimes \F_n$.

It the action is free, $A \subset M$ is the unique Cartan subalgebra up to unitary conjugacy. If the action is not free, $M$ has no Cartan subalgebras.
\end{theorem}

The proof of Theorem \ref{thm.uniqueCartan} consists of two parts. The group $\F_n$ has the \emph{complete metric approximation property} (CMAP) \cite{Ha78}. This means that there exists a sequence $\vphi_k : \Gamma \recht \C$ of finitely supported functions converging to $1$ pointwise and such that the maps $\theta_k : u_g \mapsto \vphi_k(g) u_g$ are \emph{completely bounded} with $\limsup \cbnorm{\theta_k} = 1$. Then, $\theta_k$ automatically extends to an ultraweakly continuous map $\rL \F_n \recht \rL \F_n$ without increasing $\cbnorm{\theta_k}$. Since $\F_n \actson (X,\mu)$ is profinite, it follows that $M$ has CMAP as a II$_1$ factor: there exists a sequence of finite rank, ultraweakly continuous, completely bounded maps $\theta_k : M \recht M$ converging pointwise in $\|\, \cdot \,\|_2$ to the identity and satisfying $\limsup \cbnorm{\theta_k} = 1$. Ozawa and Popa prove the following general statement: let $M$ be a II$_1$ factor with CMAP and $P \subset M$ a diffuse amenable (in particular, abelian) subalgebra; denote by $\cG$ the group of unitaries in $M$ normalizing $P$. Then, the action of $\cG$ on $P$ by automorphisms $\Ad u$ is necessarily \emph{weakly compact.} We do not explain this notion here but just note that the slightly stronger notion of \emph{compactness} means that the closure of $\Ad \cG$ inside the Polish group $\Aut P$ is compact.

So far, the reasoning in the previous paragraph works for any profinite action of a CMAP group. In the second part of the proof, Ozawa and Popa prove the following. Let $\F_n \actson (X,\mu)$ be any ergodic p.m.p.\ action, $P \subset \rL^\infty(X) \rtimes \F_n$ a subalgebra and $\cG$ a group of unitaries normalizing $P$ such that the action $\Ad \cG$ on $P$ is weakly compact. Then either $\cG$ generates an amenable von Neumann algebra or there almost exists a unitary conjugacy of $P$ into $\rL^\infty(X)$ (in the sense of Theorem \ref{thm.intertwine}). The two parts together yield Theorem \ref{thm.uniqueCartan}.

\begin{problem}
Let $n \geq 2$ and let $\F_n \actson (X,\mu)$ be \emph{any} free ergodic p.m.p.\ action. Does $\rL^\infty(X) \rtimes \F_n$ always have a unique Cartan subalgebra up to unitary conjugacy?

Let $n \geq 2$ and let $M$ be an arbitrary II$_1$ factor. Is it always true that $M \ovt \rL \F_n$ has no Cartan subalgebra? By \cite{OP07}, this is indeed the case if $M$ has the complete metric approximation property.
\end{problem}

Another breakthrough unique Cartan decomposition theorem was obtained in \cite{PV09}. On the one hand it is weaker than Theorem \ref{thm.uniqueCartan} since we are only able to deal with group measure space Cartan subalgebras, but on the other hand it is stronger since we consider arbitrary group actions.

\begin{theorem}[Popa, Vaes \cite{PV09}]\label{thm.PVCartan}
Let $\Gamma = \Gamma_1 * \Gamma_2$ be a free product where $\Gamma_1$ admits a non-amenable subgroup with the relative property (T) and where $\Gamma_2$ is any non-trivial group. Let $\Gamma \actson (X,\mu)$ be any ergodic p.m.p.\ action and put $M = \rL^\infty(X) \rtimes \Gamma$.

If the action is free, $\rL^\infty(X)$ is, up to unitary conjugacy, the unique group measure space Cartan subalgebra. If the action is not free, $M$ has no group measure space Cartan subalgebra.
\end{theorem}

Write $A = \rL^\infty(X)$. So, $M = A \rtimes \Gamma$ and we denote by $x = \sum_{g \in \Gamma} x_g u_g$ the unique Fourier expansion of $x \in M$, with $x_g \in A$ for all $g \in \Gamma$.

Assume that $M = B \rtimes \Lambda$ is any other group measure decomposition. Denote by $(v_s)_{s \in \Lambda}$ the canonical unitary elements that correspond to this decomposition. The first step of the proof of Theorem \ref{thm.PVCartan} consists in transferring some of the rigidity of $\Gamma$ to $\Lambda$. Assume that $\vphi_n$ is a deformation of the identity of $M$. If instead of $\Gamma$, our unknown group $\Lambda$ would have a non-amenable subgroup $\Lambda_0$ with the relative property (T), this would imply that on the unitary elements $v_s, s \in \Lambda_0$, the deformation $\vphi_n$ converges uniformly to the identity: for $n$ large enough and all $s \in \Lambda_0$, we get that $\|\vphi_n(v_s) - v_s\|_2$ is small. Moreover, the abelian algebra $A$ cannot contain a copy of the non-amenable algebra $\rL \Lambda_0$ so that Theorem \ref{thm.intertwine} provides a sequence $s_k$ in $\Lambda_0$ such that the Fourier coefficients of $v_{s_k}$ tend to zero pointwise in $\|\, \cdot \,\|_2$~: for all $g \in \Gamma$, we have $\|(v_{s_k})_g\|_2 \recht 0$ as $k \recht \infty$.

Obviously, we cannot prove that $\Lambda$ automatically has a non-amenable subgroup with the relative property (T). Nevertheless, we have the following transfer of rigidity result.

\begin{lemma}\label{lem.transfer}
Let $A \rtimes \Gamma = M = B \rtimes \Lambda$ be two crossed product decompositions of the same II$_1$ factor with $A$ and $B$ being amenable. Assume that $\Gamma$ admits a non-amenable subgroup with the relative property (T) and let $\vphi_n$ be a deformation of the identity of $M$. For every $\eps > 0$ there exists an $n \in \N$ and a sequence of group elements $s_k \in \Lambda$ such that
\begin{itemize}
\item $\|\vphi_n(v_{s_k}) - v_{s_k}\|_2 < \eps$ for all $k$,
\item the Fourier coefficients of $v_{s_k}$ tend to zero pointwise: for all $g \in \Gamma$, we have $\|(v_{s_k})_g\|_2 \recht 0$ as $k \recht \infty$.
\end{itemize}
\end{lemma}

Lemma \ref{lem.transfer} is proven by playing the following positive-definite ping-pong: the formula $\psi_n(s) =  \tau(v_s^* \vphi_n(v_s))$ defines a sequence of positive definite functions on $\Lambda$, which in turn define completely positive maps $\theta_n : B \rtimes \Lambda \recht B \rtimes \Lambda : \theta_n(b v_s) = \psi_n(s) b v_s$, which in their turn define positive definite functions $\gamma_n : \Gamma \recht \C : \gamma_n(g) = \tau(u_g^* \theta_n(u_g))$. By construction, $\gamma_n \recht 1$ pointwise and hence uniformly on the subgroup of $\Gamma$ with the relative property (T). From this, one can deduce the conclusion of Lemma \ref{lem.transfer}.

The starting point to prove Theorem \ref{thm.PVCartan} is an application of Lemma \ref{lem.transfer} to the word length deformation $\vphi_\rho$ on $M = A \rtimes (\Gamma_1 * \Gamma_2)$ given in Example \ref{ex.Haag-freeprod}. The fact that $\|\vphi_\rho(x) - x\|_2$ is small means that $x$ lies close to a linear combination of $a u_g$, $a \in A$ and $g \in \Gamma$ with $|g|$ not too large. We refer to such elements $x \in M$ as being \lq short\rq. Lemma \ref{lem.transfer} provides a sequence of short unitaries $v_{s_k}$. Since $B$ is abelian and normalized by the short unitaries $v_{s_k}$, a combinatorial argument implies that the elements of $B$ are themselves \lq uniformly short\rq. But Ioana, Peterson and Popa, in their study of rigid subalgebras of amalgamated free products \cite{IPP05}, proved that this implies that $B$ can be unitarily conjugated into one of the \lq obvious\rq\ short subalgebras of $A \rtimes (\Gamma_1 * \Gamma_2)$, namely $A \rtimes \Gamma_1$ or $A \rtimes \Gamma_2$. It follows that $B$ can actually be conjugated into $A$ since otherwise the normalizer of $B$, i.e.\ the whole of $M$, would get conjugated in $A \rtimes \Gamma_i$ as well.

\begin{problem}
Does the transfer of rigidity lemma hold for arbitrary Cartan subalgebras? More precisely, let $B \subset M$ be a Cartan subalgebra. Does there exist an $n$ and a sequence of unitaries $v_k \in M$, normalizing $B$ and satisfying the same two properties as the unitaries $v_{s_k}$ in the formulation of Lemma \ref{lem.transfer}?

In the affirmative case, one can replace \lq group measure space Cartan\rq\ by \lq Cartan\rq\ throughout the formulation of Theorem \ref{thm.PVCartan}.
\end{problem}

\section{Superrigidity for group measure space factors} \label{sec.superrigid}

As explained in paragraph \ref{par.superrigid}, W$^*$-superrigidity of an action $\Gamma \actson (X,\mu)$ arises as the \lq sum\rq\ of OE superrigidity and the uniqueness of the group measure space Cartan subalgebra in $\rL^\infty(X) \rtimes \Gamma$, up to unitary conjugacy. This makes W$^*$-superrigidity theorems extremely hard to obtain.

Unfortunately, none of the profinite actions covered by the uniqueness of Cartan theorems in \cite{OP07,OP08} is known to be (virtually) OE superrigid. Also actions $\Gamma_1 * \Gamma_2 \actson (X,\mu)$ are not OE superrigid so that we cannot directly apply Theorem \ref{thm.PVCartan}. But Theorem \ref{thm.PVCartan} can be generalized so that it covers arbitrary actions of certain \emph{amalgamated free products} $\Gamma = \Gamma_1 *_\Sigma \Gamma_2$ over a common amenable subgroup $\Sigma$ which is sufficiently non-normal inside $\Gamma$ (see \cite{PV09}).
In combination with Popa's OE superrigidity for Bernoulli actions \cite{Po05} or Kida's OE superrigidity \cite{Ki09}, we obtained the following examples of W$^*$-superrigid actions. A general statement and more examples, including Gaussian actions and certain co-induced actions, can be found in \cite{PV09}.

\begin{theorem}[Popa, Vaes \cite{PV09}]\label{thm.PVsuperrigid}
Let $n \geq 3$ and denote by $T_n < \PSL(n,\Z)$ the subgroup of upper triangular matrices.

Let $\Lambda \neq \{e\}$ be an arbitrary group and $\Sigma < T_n$ an infinite subgroup. Put $\Gamma = \PSL(n,\Z) *_\Sigma (\Sigma \times \Lambda)$. Then, the Bernoulli action $\Gamma \actson (X_0,\mu_0)^\Gamma$ is W$^*$-superrigid. More generally, whenever $\Gamma \actson I$ and $\Sigma \cdot i$ is infinite for all $i \in I$, the generalized Bernoulli action $\Gamma \actson (X_0,\mu_0)^I$ is W$^*$-superrigid.

Any free mixing p.m.p.\ action of $\PSL(n,\Z) *_{T_n} \PSL(n,\Z)$ is W$^*$-superrigid.
\end{theorem}

The following beautiful result was obtained very recently by Ioana \cite{Io10} and should be compared with Theorem \ref{thm.strongrigid}.

\begin{theorem}[Ioana \cite{Io10}]\label{thm.Io-superrigid}
Let $\Gamma$ be an icc group that admits an infinite normal subgroup with the relative property (T). Then, the Bernoulli action $\Gamma \actson (X_0,\mu_0)^\Gamma$ is W$^*$-superrigid.
\end{theorem}

Again, because of Popa's OE superrigidity for Bernoulli actions \cite{Po05}, the issue is to prove that $M$ has a unique group measure space Cartan subalgebra up to unitary conjugacy. Write $A = \rL^\infty\bigl(X_0^\Gamma\bigr)$ and $M = A \rtimes \Gamma$. Assume that $M = B \rtimes \Lambda$ is another group measure space decomposition. The first step of the proof of Theorem \ref{thm.Io-superrigid} is similar to the transfer of rigidity ping-pong technique that I explained after Lemma \ref{lem.transfer}. Denote by $(v_s)_{s \in \Lambda}$ the canonical group of unitaries in $B \rtimes \Lambda$. So, Ioana considers the $*$-homomorphism $\Delta : B \rtimes \Lambda \recht (B \rtimes \Lambda) \ovt \rL \Lambda$ given by $\Delta(b v_s) = b v_s \ot v_s$ for all $b \in B$, $s \in \Lambda$. Since $\rL \Lambda \subset M$, we rather view $\Delta$ as an embedding of $M$ into $M \ovt M$.

But $M$ arises from a Bernoulli action of a (relative) property (T) group. As a corollary of Theorem \ref{thm.strongrigid}, every automorphism of $M$ is of a special form, i.e.\ the composition of an inner automorphism and automorphisms induced by a character $\om : \Gamma \recht S^1$ and by a self-conjugacy of the action $\Gamma \actson A$. In an amazing technical tour de force Ioana manages to generalize Popa's methods from automorphisms to embeddings and to give an almost complete picture of all possible embeddings of $M$ into $M \ovt M$, when $M$ arises from a Bernoulli action of a property (T) group. This picture is sufficiently precise so that applied to the embedding $\Delta$ constructed in the previous paragraph, one can conclude that $A$ and $B$ are unitarily conjugate.

\end{document}